\documentclass[openany, amssymb, psamsfonts]{amsart}
\usepackage{amssymb}
\usepackage{mathrsfs,comment}
\usepackage[usenames,dvipsnames]{color}
\usepackage{amsmath}
\usepackage[normalem]{ulem}
\usepackage{url}
\usepackage[all,arc,2cell]{xy}
\UseAllTwocells
\usepackage{mathtools}
\usepackage{enumerate}
\usepackage{xcolor}
\usepackage{tabularx}
\usepackage{hyperref}

\newcommand{\cC}{\mathcal{C}}

\newcommand{\cO}{\mathcal{O}}

\newcommand{\cM}{\mathcal{M}}

\newcommand{\N}{\mathbb{N}}

\DeclareSymbolFont{extraup}{U}{zavm}{m}{n}
\DeclareMathSymbol{\varheart}{\mathalpha}{extraup}{86}
\DeclareMathSymbol{\vardiamond}{\mathalpha}{extraup}{87}
\hypersetup{%
  bookmarksnumbered=true,%
  bookmarks=true,%
  colorlinks=true,%
  linkcolor=blue,%
  citecolor=blue,%
  filecolor=blue,%
  menucolor=blue,%
  pagecolor=blue,%
  urlcolor=blue,%
  pdfnewwindow=true,%
  pdfstartview=FitBH}

\def\makeautorefname#1#2{\expandafter\def\csname#1autorefname\endcsname{#2}}
\def\equationautorefname~#1\null{(#1)\null}
\makeautorefname{footnote}{footnote}%
\makeautorefname{item}{item}%
\makeautorefname{figure}{Figure}%
\makeautorefname{table}{Table}%
\makeautorefname{part}{Part}%
\makeautorefname{answer}{Answer}%
\makeautorefname{appendix}{Appendix}%
\makeautorefname{chapter}{Chapter}%
\makeautorefname{fact}{Fact}%
\makeautorefname{section}{Section}%
\makeautorefname{subsection}{Section}%
\makeautorefname{subsubsection}{Section}%
\makeautorefname{theorem}{Theorem}%
\makeautorefname{thm}{Theorem}%
\makeautorefname{cor}{Corollary}%
\makeautorefname{lem}{Lemma}%
\makeautorefname{prop}{Proposition}%
\makeautorefname{pro}{Property}
\makeautorefname{conj}{Conjecture}%
\makeautorefname{defn}{Definition}%
\makeautorefname{notn}{Notation}
\makeautorefname{notns}{Notations}
\makeautorefname{rem}{Remark}%
\makeautorefname{quest}{Question}%
\makeautorefname{exmp}{Example}%
\makeautorefname{ax}{Axiom}%
\makeautorefname{claim}{Claim}%
\makeautorefname{ass}{Assumption}%
\makeautorefname{asss}{Assumptions}%
\makeautorefname{con}{Construction}%
\makeautorefname{prob}{Problem}%
\makeautorefname{warn}{Warning}%
\makeautorefname{obs}{Observation}%
\makeautorefname{conv}{Convention}%

\newtheorem{thm}{Theorem}[section]
\newtheorem{cor}{Corollary}[section]
\newtheorem{prop}{Proposition}[section]
\newtheorem{lem}{Lemma}[section]

\theoremstyle{definition}
\newtheorem{defn}{Definition}[section]

\newtheorem{exmp}{Example}[section]

\makeatletter
\let\c@obs=\c@thm
\let\c@cor=\c@thm
\let\c@prop=\c@thm
\let\c@lem=\c@thm
\let\c@prob=\c@thm
\let\c@con=\c@thm
\let\c@conj=\c@thm
\let\c@defn=\c@thm
\let\c@notn=\c@thm
\let\c@notns=\c@thm
\let\c@exmp=\c@thm
\let\c@fact=\c@thm
\let\c@ax=\c@thm
\let\c@pro=\c@thm
\let\c@ass=\c@thm
\let\c@warn=\c@thm
\let\c@rem=\c@thm
\let\c@sch=\c@thm
\let\c@conv=\c@thm
\let\c@equation\c@thm
\numberwithin{equation}{section}
\makeatother

\title{An Equivariant Generalization of McDuff's Theorem}

\author{Zhenghui(Sunny) Zhang}

\begin{document}
\begin{abstract}
In 1976, Kan and Thurston proved the theorem that any path-connected space $X$ is homology equivalent to the classifying space of some discrete group $G$. In 1979, McDuff proved a homotopy version of it: any path-connected space $X$ has the same weak homotopy type as the classifying space of some discrete monoid $M$. In 1984, Fiedorowicz reproved McDuff's theorem using a largely categorical construction. In this paper we will generalize Fiedorowicz's proof of McDuff's theorem to the equivariant case. Precisely, we will prove that any $G$-connected space $X$ with a $G$-fixed basepoint $x_0$ has the same weak homotopy type as the classifying space of some discrete $G$-monoid. 
\end{abstract}
\maketitle

\tableofcontents

\section{A Statement of The Main Result}
Here we state the original version of McDuff's theorem. 
\begin{thm}
Every path-connected space has the same weak homotopy type as the classifying space $BM$ of some discrete monoid $M$.
\end{thm}
We generalized this into the following equivariant case:
\begin{thm} Let $G$ be a discrete group. Every $G$-connected space $X$ with a $G$-fixed basepoint $x_0$ has the same weak homotopy type as the classifying space $BM$ of some discrete $G$-monoid $M$. Moreover, for any subgroup $H \leq G$, $X^H$ has the same weak homotopy type as $B(M^H)$, where $\cdot ^H$ denotes the fixed point set. 
\end{thm}
Toward this theorem, we will first introduce a construction of classifying spaces in Section 2. Then we will sketch the proofs by McDuff and Fiedorowicz and compare them in Section 3. Finally we will prove the equivariant version of the theorem in Section 4.
\section{A Functorial Construction of Classifying Spaces} 
In this section, we lay out a categorical construction of classifying spaces using the language of simplicial sets. 
\begin{defn}\label{defn1.1}
The category of \emph{finite ordered sets} $\Delta$ has objects as totally ordered sets $[n]=\{0, \cdots, n\}$ for any $n \in \mathbb{N}$ and morphisms as order-preserving maps $f: [m]\rightarrow [n]$. All maps in $\Delta$ are generated by compositions of \emph{face maps} $\delta_i:[n-1]\rightarrow [n]$ and \emph{degeneracy maps} $\sigma_i:[n+1]\rightarrow [n]$ defined as:
\begin{align*}
\delta_i(j)=
\begin{cases}
j &\text{ if }j<i\\
j+1 &\text{ if }j\geq i
\end{cases}\\
\sigma_i(j)=
\begin{cases}
j&\text{ if }j\leq i\\
j-1&\text{ if }j>i
\end{cases}
\end{align*}
for all $0\leq i\leq n$. In words, $\delta_i$ skips the $i^{th}$ index in the codomain; $\sigma_i$ repeats the $i^{th}$ index in the codomain. 
\end{defn}



\begin{defn}\label{defn1.14}
A \emph{simplicial object} in a category $\cC$ is a functor $\Delta^{op}\rightarrow \cC$. Classical examples are simplicial sets, simplicial spaces, and simplicial categories. 
\end{defn}

\begin{defn}\label{const}
A \emph{constant simplicial set} $X$ is a simplicial set with $X_i=X$ for all $ i \in \mathbb{N}$. All face and degeneracy maps are identities. 
\end{defn}
\begin{defn}\label{defn2.5}
A \emph{bisimplicial set} is a functor $\Delta^{op}\times\Delta^{op}\rightarrow Set$. A \emph{bisimplicial space} is a functor $\Delta^{op}\times\Delta^{op}\rightarrow Top$. 
\end{defn}
\begin{defn}\label{nerve}
 The \emph{nerve functor} $N: Cat \rightarrow sSet$ takes $N_0\mathcal{C}=\{ Ob(C)\}$. For $n\geq 1$, $N_n \mathcal{C}$ is defined as the set of all $n$-composable morphisms $(f_1, \cdots, f_n)$. The degeneracy and face maps are:
\begin{align*}
d_i(f_1, \cdots, f_n)&=
\begin{cases}
(f_2, \cdots, f_n), &i=0 \\
(f_1,\cdots,f_{i+1}\circ f_{i} ,\cdots,  f_n), & 1\leq i\leq n-1\\
(f_1, \cdots, f_{n-1}), &i=n
\end{cases}\\
s_i(f_1, \cdots, f_n)&=(f_1, \cdots, f_{i-1}, id, f_{i},\cdots, f_n)
\end{align*}
\end{defn}

\begin{defn}\label{stdsimplex}
Define the \emph{standard geometric n-simplex }$\Delta[n]^t$ as the subspace 
\[\{(t_0, ..., t_n):0\leq t_i\leq 1\text{ and }\sum_{i=0}^nt_i=1\}\] of $\mathbb{R}^{n+1}$. $\Delta[\cdot]^t$ is a functor $\Delta\rightarrow Top$. Its \emph{face maps} $\bar{\delta}_i:\Delta[n-1]^t\rightarrow \Delta[n]^t$ and \emph{degeneracy maps} $\bar{\sigma}_i:\Delta[n+1]^t\rightarrow \Delta[n]^t$ are defined as:
\begin{align*}
\bar{\delta}_i(t_0, ..., t_{n-1})&=(t_0, ..., t_{i-1}, 0, t_{i}, ..., t_{n-1}) \\ \bar{\sigma}_i(t_0, ..., t_{n+1})&=(t_0, ..., t_{i-1}, t_{i}+t_{i+1}, t_{i+2}, ..., t_{n+1})   
\end{align*}
\end{defn}

\begin{defn} \label{defn2.6}
Define the \emph{geometric realization functor} $|\cdot|:\textit{sSet}\rightarrow \textit{Top}$ as follows. For a simplicial set $X$, regard each set $X_n$ as a space with discrete topology. The topological space associated to $X$ is 
 \[|X|=\coprod_{n\geq 0}X_n\times \Delta[n]^t/\sim\]
 where the equivalence relation is $(f^*x, u)\sim (x, \bar{f}u)$ for $x \in X_n$, $u \in \Delta[m]^t$, $f:[m]\rightarrow [n]$, $f^*:X_n\rightarrow X_m$, and $\bar{f}: \Delta [m]^t  \rightarrow \Delta [n]^t $, where $f^*, \bar{f}$ are the induced maps in $X$ and $\Delta[\cdot]^t$. 
\end{defn}
\begin{defn}\label{defn2.7}
Similarly, we can define the \emph{geometric realization of simplicial spaces} by replacing the discrete set $X_n$ with space $X_n$ in Definition \ref{defn2.6}. 
\end{defn}
\begin{defn}\label{defn2.9}
Given a small category $\cC$, the classifying space $B\cC$ is defined to be $|N\cC|$, i.e. the space obtained after first applying the nerve functor then the geometric realization functor. We call $B:Cat \rightarrow Top$ the \emph{classifying space functor}. 
\end{defn}
\begin{exmp}
By regarding a group $G$ or a monoid $M$ as a category with only one object, we can define its classifying space using Definition \ref{defn2.9}. 
\end{exmp}

\section{McDuff's Theorem}
\subsection{McDuff's Proof}
In this section we sketch McDuff's original proof. The main idea is that for any simplicial complex $P$, there is an explicit way to construct a monoid $M$ based on the simplicies of $P$ such that $P \simeq BM$. 
\begin{defn}
A (discrete) \emph{monoid} is a set with an associative and unital binary operation. We use $\mathcal{M}$ to denote the category of (discrete) monoids. A monoid $M$ can also be regarded as a category with one object. 
\end{defn}
\begin{defn}\label{defn3.2}
A \emph{semigroupoid} is a nonempty, small, topologically discrete category where all the objects are isomorphic. For each semigroupoid $M_S$ and any $x \in Obj(M_S)$, there exists an associated monoid $M=Mor(x, x)$. $M$ is independent of the choice of $x$ since all elements are isomorphic in $M_S$. The category $M_S$ is isomorphic to $M\times J$, where $J$ has its set of objects as that of $M_S$, and exactly one morphism between any pair of objects as mentioned in \protect\cite{Mc}. Note that the definition of this semigroupoid is from \protect\cite{Mc} and different from the usual definition of semigroupoids. 
\end{defn}
\begin{lem}
Given a semigroupoid $M_S$ and its associated monoid $M$, the natural inclusion map $M\hookrightarrow M_S$ induces a homotopy equivalence $BM \rightarrow BM_S$. 
\end{lem}
\begin{proof}
As mentioned in Definition \ref{defn3.2}, $M_S \cong M\times J$, and since the classifying space functor preserves finite products, $BM_{S} \cong BM \times BJ$. Note that $BJ$ is contractible since $J$ has every element a zero object. This implies $BM_S\cong BM \times BJ \simeq BM$. 
\end{proof}
\begin{thm}\label{thm3.5}
For any ordered connected simplicial complex $P$, there exists a semigroupoid $F(P)$ such that $P\simeq BF(P)$. 
\end{thm}
\begin{proof}
There is a canonical way to construct $F(P)$ based on the simplicies in $P$. Let the objects in $F(P)$ correspond to the vertices in $P$. The morphisms in $F(P)$ can be generated by simplicies in $P$ with additional conditions. A precise description of this can be found in \protect\cite{Mc}. The proof of $P \simeq BF(P)$ is inductive. \\
\textit{Base Case.} If $dim(P)=0$, then $F(P)$ is the category with one object and one morphism. Suppose $dim(P)\geq 1$. For any 1-dim subcomplex $Q'\subset P$, $F(Q')$ is generated by one invertible morphism $v_0\rightarrow v_1$ for each 1-simplex $(v_0, v_1) \in Q'$, $v_0<v_1$. Therefore the canonical inclusion map $\lambda_{Q'}:Q'\rightarrow BF(Q')$ is a homotopy equivalence.\\
\textit{Induction.} Now suppose there exists a subcomplex $Q \subset P$ such that $\lambda_Q:Q\rightarrow BF(Q)$ is a homotopy equivalence. Let $K=Q\cup_{\partial \sigma}\sigma$ where $\sigma$ is an $n$-simplex with $n \geq 2$. The inclusion $Q \hookrightarrow K$ induces an inclusion $F(Q)\hookrightarrow F(K)$, and we can construct an equivalence $\lambda_K:K\rightarrow BF(K)$ such that the following diagram commutes:
\[
\xymatrix{
Q \ar[r]^-{\lambda_Q}\ar[d]_{i}& BF(Q)\ar[d]^{i_*}\\
K \ar[r]^-{\lambda_K}& BF(K)
}
\]
By induction this proves $P \simeq BF(P)$ for any finite complex $P$. For any infinite complex, using the fact that for any directed set of ordered, connected complexes $Q^\alpha$, $F(lim_{\rightarrow}Q^\alpha)=lim_{\rightarrow}F(Q^\alpha)$, this also holds. 
\end{proof}
\begin{lem}\label{lem3.5}
Every space is weak homotopy equivalent to a simplicial complex. 
\end{lem}
\begin{proof}
By the CW approximation theorem, every space is weak homotopy equivalent to a CW complex. From 2C.5 in \protect\cite{Hatcher}, every CW complex is homotopy equivalent to a simplicial complex. 
\end{proof}
\begin{cor}
Any path-connected topological space has the weak homotopy type of the classifying space of a discrete monoid. 
\end{cor}
\begin{proof}
This follows from Theorem \ref{thm3.5} and Lemma \ref{lem3.5}. 
\end{proof}
The reason that this proof does not obviously generalize to the equivariant case is because the homotopy equivalence between a simplicial complex $P$ and $BF(P)$ is proved inductively. In each inductive step there is no natural way to extend the $G$-action.  
\subsection{Fiedorowicz's Proof}
Fiedorowicz's proof built on the fact that every connected space is weak homotopy equivalent to the classifying space of a \emph{topologial monoid}. His proof shows that for any topological monoid $M$, there is a systematic way to construct another topological monoid $N$ with underlying discrete monoid $N^\delta$ such that $BM$ is weak homotopy equivalent to $BN^\delta$. Since the proof is quite involved, we give an overview here. 
\begin{enumerate}
    \item Lemma \ref{lem3.6} shows that any path-connected space $X$ is weak homotopy equivalent to the classifying space of a topological monoid $M$
    \item For each topological monoid $M$, Lemma \ref{lem3.17} constructs a discrete simplicial category $M_*^\delta$ such that there exists a weak homotopy equivalence $BM_*^\delta \rightarrow BM$. The proof of this heavily relies on Lemma \ref{lem3.14} and \ref{lem3.16}. Then the only thing left to show is a weak homotopy equivalence $BM^\delta\rightarrow BM^\delta_*$. 
    \item  Proposition \ref{lem3.18} shows it's sufficient to construct a functor $K: M^I\rightarrow M$ satisfying certain conditions to show a weak homotopy equivalence $BM^\delta\rightarrow BM^\delta_*$. This proof uses the fact that after applying the classifying space functor $B$, natural transformations become homotopies
    \item Apply (3) to the specific setting. Given a weak homotopy equivalence $X \rightarrow BM$, Theorem \ref{thm3.9} constructs a topological monoid $N$ homotopic to $M$ and a functor $K:N^{I} \rightarrow N$ satisfying all conditions. Then Theorem \ref{thm3.12} shows weak homotopy equivalences $BN^\delta \rightarrow BN \leftarrow BM \rightarrow X$
\end{enumerate}
In the following we provide necessary definitions and proofs of the claims. 
\begin{defn}
A \emph{topological monoid} $M$ is a monoid with a topology such that the multiplication $M \times M\rightarrow M$ is continuous. The \emph{basepoint} of a topological monoid is the identity element. We use $\mathcal{TM}$ to denote the category of topological monoids. Denote $\cdot ^\delta: \mathcal{TM}\rightarrow \mathcal{M}$ as the forgetful functor that forgets the topology. 
\end{defn}
\begin{defn}
Given two topological spaces $X, Y$, denote $Y^X$ as the space $Map(X, Y)$ with compact open topology. 
\end{defn}
\begin{defn}
Given a pointed space $(X, x_0)$, the basepoint is \emph{non-degenerate} if the inclusion map $x_0 \hookrightarrow X$ is a cofibration. 
\end{defn}
\begin{defn}\label{moore}
The \emph{Moore loop space} $\Lambda X$ on $X$ is a subspace of $X^{[0, \infty)}$. Given a fixed basepoint $x_0 \in X$, its elements are pairs $(\alpha, r)$ where $\alpha: [0, \infty)\rightarrow X$, $\alpha(t)=x_0$ for $t=0$ and $t\geq r$. The space $\Lambda X$ can be regarded as a monoid by forgetting the length $r$ of each element. Define the multiplication of elements $f, g \in \Lambda X$ by the following concatenation: if $r_1$ is the value after which $f$ remains constant and $r_2$ is the value for $g$, then $g\circ f$ is:
\[g\circ f(x)=\begin{cases}
f(x) & \text{for $0\leq x\leq r_1$}\\
g(x-r_1)& \text{for $r_1\leq x\leq r_1+r_2$}\\
x_0 & x\geq r_1+r_2
\end{cases}\]
This multiplication structure makes $\Lambda X$ a topological monoid. 
\end{defn}
\begin{lem}\label{lem3.6}
Any connected space $X$ is weak homotopy equivalent to the classifying space of the topological monoid $M=\Lambda X$. 
\end{lem}
\begin{proof}
By the CW approximation theorem, we can replace the space $X$ with a CW complex of the same weak homotopy type functorially. Choose the basepoint of $X$ as a $0$-cell of $X$. Then Lemma 15.4 in \protect\cite{May2} shows that the map $\xi:B \Lambda X \rightarrow X$ defined by 
\[\xi|[\lambda_1, \cdots, \lambda_p], u|=(\lambda_1\cdots \lambda_p)(\sum_{i=1}^p u_il(\lambda_i))\]
for $\lambda_i \in \Lambda X$, $u=(t_0, \cdots, t_p) \in \Delta[p]$, and $u_i=t_0+\cdots+t_{i-1}$ induces a weak homotopy equivalence. Note that because $\Lambda X$ is a topological monoid, the nerve functor sends $\Lambda X$ to a simplicial space and the geometric realization $|N\Lambda X|=B\Lambda X$ is that of the simplicial space as in Definition \ref{defn2.7}. 
\end{proof}
\begin{defn}
Given two topological monoids $M, N$ where $N$ has an action on the space $X$, define the \emph{wreath product} $N \int M^X$ by taking the topology as the product topology on $N \times M^X$. For the monoid structure, given two elements $(\alpha, f)$, $(\beta, g)$,  
\[(\alpha, f)(\beta, g)=(\alpha \beta, u)\]
where $u(x)=f(\beta x)\cdot g(x)$. Hence $N \int M^X$ is a topological monoid. 
\end{defn}
\begin{defn}
Let $H(X)$ denote the topological monoid of self homotopy equivalences of the space $X$. 
\end{defn}
\begin{defn}
A category \emph{enriched} over topological spaces is a \emph{topological category}. 
\end{defn}
\begin{defn}\label{const}
Given a topological category $\cC$ and a topological space $X$, $\cC^X$ is a \emph{topological function category} with objects $\{X\rightarrow Ob(\cC)\}$ and morphisms $\{X\rightarrow Mor(\cC)\}$. We denote 
\[J: \cC \rightarrow \cC^X\]
as the constant functor. 
\end{defn}
\begin{lem}\label{fact3.10}
Let $\pi: X \times I\rightarrow X$ be the projection map and $\phi: H(X)\rightarrow H(X\times I)$ be the natural inclusion map $\alpha \mapsto \alpha \times id$. Then there is an induced monoid homomorphism:
\[\phi \int \pi^*: H(X)\int M^X \rightarrow H(X\times I)\int M^{X\times I}\]
\end{lem}
\begin{lem}\label{lem3.14}
Let $f_*: X_* \rightarrow Y_*$ be a map of (Reedy cofibrant) simplicial spaces such that each $f_n:X_n\rightarrow Y_n$ is a (weak) homotopy equivalence. Then $|f_*|:|X_*|\rightarrow |Y_*|$ is a (weak) homotopy equivalence. 
\end{lem}
The Reedy cofibrancy is a mild cofibration condition that will be satisfied in our context. 
\begin{lem}\label{lem3.15}
If $X_{**}$ is a bisimplicial space as defined in Definition \ref{defn2.5}, there are natural homeomorphisms 
\[
|m\mapsto |n\mapsto X_{mn}||\cong
|m\mapsto X_{m, m}|\cong 
|n\mapsto |m\mapsto X_{mn}||
\]
That is, the order of geometric realizations doesn't affect the resulting space. 
\end{lem}
\begin{lem}\label{lemconst}
Given a constant simplicial space $X$, $|X|\cong X$. 
\end{lem}
\begin{lem}\label{lem3.16}
Consider the adjunction between simplicial sets and topological spaces, where $Sing(-)$ stands for the totally singular functor and $|\cdot|$ stands for the geometric realization functor. 
\[
\xymatrix{
sSet \ar@<1ex>[r]^-{|\cdot|}& Top\ar@<1ex>[l]^-{Sing(-)}
}
\]
Then for any $X \in Top$, the counit map $\epsilon: |Sing(X)|\rightarrow X$ is a weak homotopy equivalence. A proof of this can be found in Theorem 16.6 in \protect\cite{May3}.
\end{lem}
\begin{lem}\label{lem3.17}
Let $M$ be a topological monoid with a non-degenerate basepoint. Regard $M$ as a topological category. Suppose that for every $n\in \N$,
\[J^\delta: M^\delta \rightarrow (M^{\Delta^n})^\delta\]
induces a homotopy equivalence of classifying spaces. Then the natural map $BM^\delta\rightarrow B M$ is a weak homotopy equivalence. 
\end{lem}
\begin{proof}
Define a discrete, simplicial category $M^\delta_*: \Delta^{op}\rightarrow Cat$ as $M^\delta_n=M^\delta_*([n])=({M^{\Delta^n}})^\delta$. $({M^{\Delta^n}})^\delta$ is a category as in Definition \ref{const}. Then the nerve of $M^\delta_*$ is a bisimplicial set $N_*(M^\delta_*)$ such that degree $(m,n)$ is:
\[N_m(M^\delta_n)=\{\Delta^n \rightarrow \underbrace{M \times \cdots \times M}_{m \text{ times}}\}^\delta\]
This is because by definition $N_m(M^\delta_n)$ is the set of $m$-tuples of composable morphisms in $M^\delta_n$ and $M^\delta_n$ is a category with only one object. By Lemma \ref{lem3.15}, 
\[BM_*^\delta=|N(M_*^\delta)|=|m \mapsto |n \mapsto N_m(M_n^\delta)|| \]
Fixing $m$, $N_m(M_*^\delta)$ is a simplicial set equal to $Sing(N_m M)$ by an easy comparison of definitions. Therefore 
\[BM_*^\delta=|m \mapsto |Sing(N_m M)|| \]
Also
\[BM=|m \mapsto N_m M| \]
By Lemma \ref{lem3.16}, $|Sing(N_m M)|$ is weak homotopy equivalent to $N_mM$. The condition of $M$ having a non-degenerate basepoint guarantees that $N_*M$ and $|Sing(N_* M)|$ are Reedy cofibrant simplicial spaces. Thus the condition of Lemma \ref{lem3.14} is satisfied, and there exists a weak homotopy equivalence $BM_*^\delta\rightarrow BM$. \\
Now we show $BM^\delta \rightarrow BM_*^\delta$ is a weak homotopy equivalence. By assumption, $J^\delta: BM^\delta \rightarrow B(M^{\Delta^n})^\delta$ is an equivalence. Again by Lemma \ref{lem3.15}, 
\[BM_*^\delta=|n \mapsto |m \mapsto N_m(M_n^\delta)||=|n \mapsto BM_n^\delta|=|n \mapsto B(M^{\Delta^n})^\delta|\]
Regard $BM^\delta$ as a constant simplicial space. Lemma \ref{lemconst} shows 
\[
BM^\delta=|BM^\delta|=|n \mapsto BM^\delta|
\]
Both $BM^\delta$ and $B({M^{\Delta}}^n)^\delta$ are Reedy cofibrant simplicial spaces given $M$ having a non-degenerate basepoint. Since $BM^\delta \rightarrow B({M^{\Delta}}^n)^\delta$ is a homotopy equivalence, Lemma \ref{lem3.14} implies $BM^\delta \rightarrow BM_*^\delta$ a homotopy equivalence.  Combining these two results together, $BM^\delta \rightarrow BM$ is weak homotopy equivalence. 
\end{proof}
\begin{prop}\label{lem3.18} Let $M$ be a topological monoid with a non-degenerate basepoint, and let $J$ be the constant functor as in Definition \ref{const}. If there exists a continuous functor $K: M^I \rightarrow M$ and a continuous natural transformation $\psi$ from
\[M^I \xrightarrow{K}M\xrightarrow{J}M^I\]
to the identity functor, then $BM^\delta$ is weak homotopy equivalent to $BM$. 
\end{prop}
\begin{proof}
After iterating $K$, there exists a natural transformation from
\[M^{I^n}\xrightarrow{K^n}M\xrightarrow{J}M^{I^n}\]
to the identity functor. Applying the forgetful functor, $\psi^\delta$ still a natural transformation from
\[(M^{I^n})^\delta\xrightarrow{K^n}M^\delta\xrightarrow{J^\delta}(M^{I^n})^\delta\]
to the identity functor. Let $I:(M^{I^n})^\delta \rightarrow M^\delta$ be any evaluation map, then 
\[M^\delta\xrightarrow{J^\delta}(M^{I^n})^\delta \xrightarrow{I}M^\delta\]
is the identity map. After applying the classifying space functor, natural transformations become homotopies. Then 
\begin{align*}
 BK^n \circ BJ^\delta \simeq id \\
BJ^\delta \circ BI=id
\end{align*}
Thus $BJ^\delta: BM^\delta \rightarrow B(M^{I^n})^\delta$ is a homotopy equivalence. By $\Delta^n\cong I^n$ and Lemma \ref{lem3.17}, there exists a weak homotopy equivalence $BM^\delta\rightarrow B M$. 
\end{proof}
\begin{thm}\label{thm3.9}
Given a path-connected space $X$, let $\overline{\Lambda X}$ denote the colimit of the following sequence where the maps are as defined in Lemma \ref{fact3.10}:
\[\Lambda X \hookrightarrow H(I)\int (\Lambda X)^I \hookrightarrow H(I^2)\int (\Lambda X)^{I^2}\hookrightarrow \cdots\]
 Then there are a pair of weak homotopy equivalences $B(\overline{\Lambda X})^\delta  \rightarrow B(\overline{\Lambda X})\leftarrow B(\Lambda X)$. 
\end{thm}
\begin{proof}
In the construction of $\Lambda X$ as in Lemma \ref{lem3.6}, we replace $X$ functorially with a CW complex with the same weak homotopy type and choose the basepoint as a $0$-cell of it. Then $X$ is a CW complex with a non-degenerate basepoint. The inclusion map $\Lambda X \rightarrow \overline{\Lambda X}$ is a homotopy equivalence, so $B(\Lambda X)\rightarrow B(\overline{\Lambda X})$ is a homotopy equivalence. Now the only thing to show is a weak homotopy equivalence $B(\overline{\Lambda X})^\delta \rightarrow B(\overline{\Lambda X})$. Since $X$ has a non-degenerate basepoint, so do $\Lambda X$ and $\overline{\Lambda X}$ by construction. Here we apply Proposition \ref{lem3.18}, which shows that it is sufficient to construct a functor 
\[K: (\overline{\Lambda X})^{I^m}\rightarrow \overline{\Lambda X}\]
that satisfies the conditions in Proposition \ref{lem3.18}. To construct this, we can instead construct a family of compatible functors $K_n:(H(I^n)\int (\Lambda X)^{I^n})^{I^m}\rightarrow H(I^{m+n})\int (\Lambda X)^{I^{m+n}}$ with respect to the sequence whose colimit is $\overline{\Lambda X}$. Then $\{K_n\}$ induces a $K: \overline{\Lambda X}^{I^m} \rightarrow \overline{\Lambda X}$. The construction of $K_n$ involves manipulating the indices in $I^{*}$ and exploits the fact that we are taking the colimit of the sequence so there is enough space to move around. A detailed construction can be found in Theorem 3.4 in \protect\cite{Fied}.  
\end{proof}
\color{black}
\begin{thm}\label{thm3.12}
Let $\mathcal{P}$ denote the category of path-connected spaces. There exists a functor $D: \mathcal{P}\rightarrow \mathcal{M}$ such that for $X \in \mathcal{P}$, $BDX$ is weak homotopy equivalent to $X$.
\end{thm}
\begin{proof}
Based on Theorem \ref{thm3.9}, there exists a pair of weak homotopy equivalences $B(\overline{\Lambda X})^\delta \rightarrow B(\overline{\Lambda X})\leftarrow B(\Lambda X)$ under the canonical map. By Lemma \ref{lem3.6}, there exists a weak homotopy equivalence $\xi: B (\Lambda X) \rightarrow X$. Since the construction of $\Lambda X$, $(\overline{\Lambda X})^\delta$ and classifying spaces are all functorial, take $DX=(\overline{\Lambda X})^\delta$. Then we have the natural weak homotopy equivalences as follows:
\[
\xymatrix{
BDX=B(\overline{\Lambda X})^\delta
\ar[r]
& B(\overline{\Lambda X})
& B \Lambda X\ar[l] \ar[r]
&  X\\
}
\]
\end{proof}
Fiedorowicz's proof of McDuff's Theorem emphasizes on categorical constructions. If we can verify that the constructions are compatible with $G$-actions and well-behaved with respect to the set of $H$-fixed points for all $H \leq G$, then there is a way to generalize it equivariantly. 
\section{An Equivariant Generalization}
\begin{defn}
A monoid $M$ with a left $G$-action by monoid homomorphisms is called a \emph{$G$-monoid}. Maps between $G$-monoids are $G$-monoid homomorphism $f$ such that for $g \in G$, $m \in M$, $f(g\cdot m)=g\cdot f(m)$. We denote this category as $\cM_G$. 
\end{defn}
\color{black}
\begin{defn}
A $G$-space is a topological space equipped with a left $G$-action that is continuous. A $G$-map between $G$-spaces $X, Y$ is a continuous map $f: X \rightarrow Y$ such that $f(g\cdot x)=g\cdot f(x)$ for all $g \in G$. A \emph{homotopy} between $G$-maps $f$, $g: X \rightarrow Y$ is a $G$-map $H: X \times I \rightarrow Y$ with $H(x, 0)=f(x)$, $H(x, 1)=g(x)$. 
\end{defn}
\begin{defn}
A topological $G$-monoid is a $G$-monoid with a topology such that the action of $G$ is continuous. 
\end{defn}
\begin{defn}Given a $G$-space $X$, let $X^H$ denote the set of fixed points of $H$, i.e. $X^H=\{x \in X | h \cdot x =x \text{ }\forall h \in H\}$. A space $X$ is called \emph{$G$-connected} if for all $H \leq G$, $X^H$ is connected. 
\end{defn}
\begin{defn}
Given a $G$-space $X$, a point $x \in X$ is called \emph{$G$-fixed} if $g\cdot x=x$ for all $g \in G$. 
\end{defn}

\begin{defn}
A $G$-map $f: X \rightarrow Y$ is a \emph{weak homotopy equivalence} if the map $f^H:X^H \rightarrow Y^H$ is a weak homotopy equivalence for all $H \leq G$. 
\end{defn}
\begin{thm}\label{cwapprox}
Given any $G$-space $X$, there exists a $G$-CW complex $\Gamma X$ such that there exists a weak homotopy equivalence $\tau:\Gamma X \rightarrow X$. Moreover $\Gamma$ is functorial. This is the CW approximation theorem in the equivariant case. 
\end{thm}
\begin{proof}
A precise proof of this can be found in Theorem 3.6 in \protect\cite{May1}. 
\end{proof}
Using this theorem, from now on we assume our $G$-space $X$ as a $G$-CW complex with basepoint $x_0$ a 0-cell. 
\begin{lem}
Given topological monoids $N$, $A$, $B$, an action of $N$ on $A$, $B$, and an $N$-map $h: A \rightarrow B$ that is also a monoid homomorphism, it induces a monoid homomorphism 
\[h_*:N \int A \rightarrow N \int B\]
\end{lem}
\begin{proof}
Consider $(n, a)$, $(n', a')\in N \int A$. We have 
\begin{align*}
h_*((n, a)\cdot (n', a'))
&=h^*((nn', (n' \cdot a)\cdot a'))\\
&=(nn', h\bigl((n' \cdot a)\cdot a'\bigl))\\
h_*(n, a)\cdot h^*(n', a')
&=(n, h(a))\cdot (n', h(a'))\\
&=(nn', \bigl(n'\cdot h(a)\bigl)\cdot h(a'))
\end{align*}
and because $h$ is a monoid homomorphism and an $N$-map, these two are the same. 
\end{proof}
\begin{cor}\label{cor4.9}
Given an inclusion map $i: M\rightarrow M'$ where $M, M'$ are $G$-monoids, $G$ acts trivially on $N$ and $X$, and $X$ is an $N$-space, there exists an induced $G$-monoid homomorphism $i_*:N \int M^X \rightarrow N \int M'^X$. 
\end{cor}
\begin{lem}\label{lem4.8}
Given $(X, x_0)$ a $G$-space with $x_0$ $G$-fixed, let $\Lambda X$ be the Moore loop space as in Definition \ref{moore}. Then $\Lambda X$, $\overline{\Lambda X}$, $(\overline{\Lambda X})^\delta$ are also $G$-spaces. 
\end{lem}
\begin{proof}
The reason we need a pointed space $X$ with $x_0$ $G$-fixed is that the construction of $\Lambda X$ requires a choice of the basepoint, and the $G$-action is well-defined only when $x_0$ is $G$-fixed. Because $\Lambda X$ is a subspace of $Map([0, \infty), X)$, given $f \in \Lambda X$, define $g\cdot f$ as the loop such that given any $t\in [0, \infty)$,
\[(g\cdot f)(t)=g\cdot f(t)\]
Recall the construction of $\overline{\Lambda X}$ in Theorem \ref{thm3.9}. We can define a $G$-action on $\overline{\Lambda X}$ as the one induced from $\Lambda X$ with trivial $G$-action on $I$. This is well-defined by Corollary \ref{cor4.9}. The $G$-action on $(\overline{\Lambda X})^\delta$ follows from that on $\overline{\Lambda X}$. 
\end{proof}
\begin{cor}
The $G$-action defined in Lemma \ref{lem4.8} makes $\Lambda X$, $\overline{\Lambda X}$, $(\overline{\Lambda X})^\delta$ topological $G$-monoids. 
\end{cor}
\begin{lem}\label{lem4.12}
Given $(X, x_0)$ a $G$-space with $x_0$ $G$-fixed, $B\Lambda X$, $B\overline{\Lambda X}$, $B(\overline{\Lambda X})^\delta$ are also $G$-spaces. 
\end{lem}
\begin{proof}
In Lemma \ref{lem4.8}, we have shown that $\Lambda X$, $\overline{\Lambda X}$ are $G$-spaces. Recall $B\Lambda X=|N\Lambda X|$, where $N$ is the nerve functor and $|\cdot|$ is the geometric realization functor. By Definition \ref{defn2.6}, we know the geometric realization of $N \Lambda X$ is
\[|N\Lambda X|=\coprod_{n\geq 0}(N\Lambda X)_n \times \Delta[n]^t/\sim\]
where the equivalence relation is $(f^*x, u)\sim (x, \bar{f}u)$ for $x \in X_n$, $u \in \Delta[m]^t$, $f:[m]\rightarrow [n]$, $f^*:(N\Lambda X)_n\rightarrow (N\Lambda X)_m$, and $\bar{f}: \Delta [m]^t  \rightarrow \Delta [n]^t $, where $f^*, \bar{f}$ are the induced maps in $N \Lambda X$ and $\Delta[\cdot]^t$. The following is basically proving that the geometric realization of a simplicial $G$-space is a $G$-space. To define the $G$-action on $|N\Lambda X|$, it is sufficient to define the $G$-action on $(N\Lambda X)_n$ as the one induced from $\Lambda X$ and take the trivial action on $\Delta[n]^t$. The only thing that is left to show is the quotient relation is compatible with the $G$-action. Given $(f^*x, u)\sim (x, \bar{f}u)$, we show the orbit of $(f^*x, u)$ equals to the orbit of $(x, \bar{f}u)$. It is a fact that $f^*$ can be written as a composition of maps $d_i, s_j$ as in Definition \ref{nerve}. Observe that $d_i, s_j$ are $G$-equivariant maps in our context. Then given $g \in G$, 
\begin{align*}
g\cdot (f^*x, u)&=( g\cdot f^*x, u)\\
&=(f^*(g\cdot x), u)\\
&=(g\cdot x, \bar{f}u)\\
&=g\cdot (x, \bar{f}u)
\end{align*}
Hence this is a well-defined $G$-action on $B\Lambda X$. Similarly we can construct the induced $G$-action on $B\overline{\Lambda X}$ and $B(\overline{\Lambda X})^\delta$. 
\end{proof}



\begin{lem}\label{lem4.13}
Recall the construction of the map $\xi: B\Lambda X \rightarrow X$ in Lemma \ref{lem3.6}. If $(X, x_0)$ is a $G$-space with $x_0$ $G$-fixed, then $\xi$ is a $G$-map. 
\end{lem}
\begin{proof}
Given an element $(|[\lambda_1, \cdots, \lambda_p], u|)$ where $\lambda_i \in \Lambda X$, $u=(t_0, \cdots, t_p) \in \Delta_p$, 
\begin{align*}
g\cdot \xi(|[\lambda_1, \cdots, \lambda_p], u|)&=g\cdot \bigl((\lambda_1 \cdots \lambda_p)\Sigma_{i=1}^p u_i l(\lambda_i)\bigl)\\
&=\bigl(g\cdot (\lambda_1 \cdots \lambda_p)\Sigma_{i=1}^p u_i l(\lambda_i)\bigl)\\
&=\bigl(g\cdot (\lambda_1 \cdots \lambda_p)\Sigma_{i=1}^p u_i l(g\cdot \lambda_i)\bigl)\\
&=\xi(g\cdot|[\lambda_1, \cdots, \lambda_p], u|)
\end{align*}
\end{proof}

\begin{thm}\label{thm4.16}
Let $\mathcal{P}_G$ denote the category of pointed $G$-connected spaces where the basepoint is $G$-fixed. There exists a functor $D: \mathcal{P}_G\rightarrow \mathcal{M}_G$ such that for $X \in \mathcal{P}_G$, $BDX$ is weak homotopy equivalent to $X$ in the equivariant sense. 
\end{thm}
\begin{proof}
Let $BDX=B(\overline{\Lambda X})^\delta$. Take the $G$-action on $B(\overline{\Lambda X})^\delta$, $B(\overline{\Lambda X})$, and $B \Lambda X$ as defined in Lemma \ref{lem4.12}. The canonical maps $B(\overline{\Lambda X})^\delta\rightarrow  B(\overline{\Lambda X})\leftarrow B\Lambda X$ are $G$-maps, and by Lemma \ref{lem4.13} $B\Lambda X \rightarrow X$ is also a $G$-map. The following $G$-maps induce weak homotopy equivalences between $BD(X^H)$ and $X^H$ by Theorem \ref{thm3.12} for all $H \leq G$. 
\[
\xymatrix{
BD(X^H)=B(\overline{\Lambda X^H})^\delta
\ar[r]
& B(\overline{\Lambda X^H})
& B \Lambda X^H\ar[l] \ar[r]
&  X^H\\
}
\]
Now we show $BD(X^H)=(BDX)^H$ for all $H \leq G$. Given any point 
\[([\{\alpha_1, \beta_1\}, \cdots, \{\alpha_k, \beta_k\}], p) \in (BDX)^H=(B(\overline{\Lambda X})^\delta)^H\]
where $[\{\alpha_1, \beta_1\}, \cdots, \{\alpha_k, \beta_k\}] \in (N\overline{\Lambda X})_k$ and $p \in \Delta[k]$. Each pair $\{\alpha_i, \beta_i\}$ denotes the two coordinates of $\overline{\Lambda X}$ and elements $\alpha_i, \beta_i$ have finite support. The point $[\{\alpha_1, \beta_1\}, \cdots, \{\alpha_k, \beta_k\}]$ is fixed under $H$ is equivalent to $img(\beta_i)\subset X^H$ for all $1\leq i \leq k$, implying $(BDX)^H=BD(X^H)$. Therefore $(BDX)^H$ is weak homotopy equivalent to $X^H$ for all $H\leq G$. What we have shown here is that all the constructions commute with fixed points. Therefore the $G$-maps induce instances of Fiedorowicz's result on passage to fixed points for each $H \leq G$. 
\end{proof}
\color{black}
In the following we show that the functor $BD(-)$ can be also understood as a natural transformation between two functors on the orbit category, each of which evaluates to $\cdot ^H$ and $BD(\cdot^H)$ on $H\leq G$. 
\begin{defn}
Define the orbit category $G\cO$ with its objects as $G$-sets $G/H$ for $H\leq G$ and morphisms as $G$-maps $\alpha: G/H \rightarrow G/K$.
\end{defn}
\begin{defn}
Given a $G$-connected space $X$ with a $G$-fixed basepoint $x_0$, define the functor $R_X: G \cO^{op} \rightarrow Top$ as $R_X(G/H)=X^H$. Given a $G$-map $\alpha: G/H \rightarrow G/K$ such that $\alpha(gH)=g \gamma K$, it induces a map 
\[R(\alpha): X^K \rightarrow X^H\]
by sending $x \mapsto \gamma^{-1} x$. 
\end{defn}
\begin{defn}
Given a $G$-connected space $X$ with a $G$-fixed basepoint $x_0$, define the functor $\bar{R}_X: G \cO^{op} \rightarrow Top$ as $\bar{R}_X(G/H)=(BDX)^H$. Given a $G$-map $\alpha: G/H \rightarrow G/K$ such that $\alpha(gH)=g \gamma K$, it induces a map 
\[\bar{R}(\alpha): (BDX)^K \rightarrow (BDX)^H\]
by sending $x \mapsto \gamma^{-1} x$. 
\end{defn}
\begin{cor}
Let $B$ denote the classifying space functor and $D$ denote the functor in Theorem \ref{thm4.16}. Then $BD:{R}_X\rightarrow \bar{R}_X$ is a natural transformation. 
\end{cor}
\begin{proof}
Given $G/H \in G\cO^{op}$, $BD(R_X(G/H))=BD(X^H)=BDX^H$. Given a $G$-map $\alpha: G/H\rightarrow G/K$ such that $\alpha(gH)=g\gamma K$, we have the following left diagram commutes as we have proved that right diagram commutes since $BD(-)$ respects the $G$-action. 
\[
\xymatrix{
R_X(G/H) \ar[d]_{BD_{R_X(G/H)}}& R_X(G/K)\ar[l]_{R_X(\alpha)}\ar[d]^{BD_{R_X(G/K)}}\\
\bar{R}_X(G/H) & \bar R_X(G/K)\ar[l]^{\bar{R}_X(\alpha)}\\
}
\Longleftrightarrow
\xymatrix{
X^H \ar[d]_{BD_{X^H}}& X^K\ar[l]_{R_X(\alpha)}\ar[d]^{BD_{X^K}}\\
(BDX)^H & (BDX)^K\ar[l]^{\bar{R}_X(\alpha)}\\
}
\]
Thus $BD(-): R_X \rightarrow \bar{R}_X$ is a natural transformation such that when evaluating on $H\leq G$, $BD_{X^H}:X^H\rightarrow (BDX)^H$. 
\end{proof}

\section*{Acknowledgments} 
I want to thank Peter May for being an extremely knowledgeable and accessible mentor who made this paper and the REU possible. I am grateful for all my friends in the REU. Thank you all for your love and support $\varheart$.  

\end{document}